\documentclass[11pt,a4paper]{article}
\usepackage{amsmath,graphicx,amsthm,amsfonts,amssymb,amstext}
\usepackage{fullpage}
\usepackage{ragged2e, enumerate,graphicx,lscape}
\usepackage{pgf, tikz, float,subcaption}
\usetikzlibrary{graphs, graphs.standard}
\usepackage{subfiles}
\usepackage{hyperref}

\newtheorem{theorem}{Theorem}

\newtheorem{lemma}{Lemma}

\newtheorem{proposition}{Proposition}

\newtheorem{conjecture}{Conjecture}

\theoremstyle{remark}

\usepackage{hyperref}
\hypersetup{
	colorlinks=true,
	linkcolor=blue,
    citecolor=purple,
    urlcolor=purple,
    }
\DeclareMathOperator{\tr}{tr}

\title{Vertex Partitioning and $p$-Energy of Graphs}
\author{Saieed Akbari \and Hitesh Kumar \and Bojan Mohar \and Shivaramakrishna Pragada}
\date{}

\begin{document}
\maketitle

\begin{abstract}
\noindent
For a Hermitian matrix $A$ of order $n$ with eigenvalues $\lambda_1(A)\ge \cdots\ge \lambda_n(A)$, define
\[ \mathcal{E}_p^+(A)=\sum_{\lambda_i > 0} \lambda_i^p(A), \quad \mathcal{E}_p^-(A)=\sum_{\lambda_i<0} |\lambda_i(A)|^p,\]
to be the positive and the negative $p$-energy of $A$, respectively. In this note, first we show that if $A=[A_{ij}]_{i,j=1}^k$, where $A_{ii}$ are square matrices, then
\[ \mathcal{E}_p^+(A)\geq \sum_{i=1}^{k} \mathcal{E}_p^+(A_{ii}), \quad \mathcal{E}_p^-(A)\geq \sum_{i=1}^{k} \mathcal{E}_p^-(A_{ii}),\]
for any real number $p\geq 1$. We then apply the previous inequalities to establish lower bounds for $p$-energy of the adjacency matrix of graphs.
\end{abstract}

\noindent
\textbf{Keywords:} Adjacency matrix, Graph energy, Positive $p$-energy, Negative $p$-energy, Schatten $p$-norm, Vertex partition

\noindent 
\textbf{MSC:} 05C50

\section{Introduction}

Let $A$ be a Hermitian matrix of order $n$. Then $A$ has $n$ real eigenvalues, which we can arrange in non-increasing order, i.e.,
\[\lambda_1(A)\ge \lambda_2(A)\ge \cdots \ge \lambda_n(A).\] 
We denote by $n^+(A)$ and $n^-(A)$ (or simply $n^+$ and $n^-$ if $A$ is clear from context) the number of positive and the number of negative eigenvalues, respectively. 
For any real number $p\ge 1$, we define 
\[\mathcal{E}_p^+(A)=\sum_{i=1}^{n^+} \lambda_i^p(A), \quad \mathcal{E}_p^-(A)=\sum_{i=n-n^-+1}^n |\lambda_i(A)|^p,\]
 and call them the \emph{positive p-energy} and the \emph{negative p-energy} of $A$, respectively. We define
 \[\mathcal{E}_p(A)=\mathcal{E}_p^+(A) + \mathcal{E}_p^-(A),\]
and call it the \emph{p-energy} of $A$. 

For a graph $G$ with adjacency matrix $A(G)$, we simply write $\mathcal{E}_p(G), \mathcal{E}_p^+(G), \mathcal{E}_p^-(G)$ to denote $\mathcal{E}_p(A(G)), \mathcal{E}^+_p(A(G)), \mathcal{E}^-_p(A(G))$, respectively. The usual \emph{graph energy} $\mathcal{E}_1(G)$ is a well-studied parameter and finds applications in many areas, including mathematical chemistry; see the survey \cite{gutman2017survey}. In recent years, the positive and negative square energies of graphs have been used to obtain good bounds for the chromatic number \cite{Ando_Lin_2015}, fractional chromatic number \cite{Guo2024new}, and vector chromatic number \cite{Coutinho2024conic}. In general, there has been a growing interest in investigating higher energies $\mathcal{E}_p(G)$ of graphs, and many interesting questions have been posed; see \cite{Elphick_FGW_2016, Nikiforov2012extremal, Nikiforov2016beyond, Tang2025p_energy_Nikiforov}. Recently, Elphick, Tang and Zhang \cite{Elphick_Tang_Zhang_2025_penergy} used $p$-energies to obtain bounds for chromatic number, which significantly generalizes many of the previously known results. 

For a Hermitian matrix $A$, the singular values of $A$ correspond to the absolute values of its eigenvalues. Consequently, the quantity $\mathcal{E}_p(A)^{1/p}$ coincides with the \emph{Schatten $p$-norm} of the matrix $A$. In \cite{Nikiforov2016beyond}, Nikiforov investigated the Schatten $p$-norms of graphs and asked which graphs are extremal with respect to the Schatten $p$-norm. Building upon this work, Tang, Liu, and Wang \cite{Tang_Liu_Wang_2025_edgeaddition} initiated the study of the positive and negative $p$-energies of graphs and posed several conjectures, one of which is stated below.

We use standard graph terminology and notation: $P_n$ denotes a path, $K_n$ denotes a complete graph, and $C_n$ denotes a cycle of order $n$. For a subset $S\subseteq V(G)$, $G[S]$ denotes the subgraph induced by $S$ in the graph $G$.

\begin{conjecture}[\cite{Tang_Liu_Wang_2025_edgeaddition}]\label{conj:E_q+}
For every connected graph $G$ of order $n$ and $p\geq 2$, 
    \[ \mathcal{E}^+_p(G)\ge \mathcal{E}^+_p(P_n).\]
\end{conjecture}

A similar conjecture can be posed for the negative $p$-energy as well.

\begin{conjecture}\label{conj:E_q-}
For every connected graph $G$ of order $n$ and $p\geq 2$, 
\[ \mathcal{E}^-_p(G)\ge \mathcal{E}^-_p(K_n).\]
\end{conjecture}

The above conjectures generalize the following seemingly innocent conjecture for connected graphs posed by 
Elphick, Farber, Goldberg and Wocjan \cite{Elphick_FGW_2016}. 

\begin{conjecture}[\cite{Elphick_FGW_2016}] \label{conj:s_plus_main}
For every connected graph $G$ of order $n$, 
    \[ \min\{\mathcal{E}_2^+(G), \mathcal{E}_2^-(G)\}\ge n-1.\]
\end{conjecture}

Only recently, linear lower bounds were obtained for $\min\{\mathcal{E}_2^+(G), \mathcal{E}_2^-(G)\}$. Zhang \cite{Zhang_2024} showed that $\min\{\mathcal{E}_2^+(G), \mathcal{E}_2^-(G)\}\ge n-\gamma(G)\ge \frac{n}{2}$, where $\gamma(G)$ denotes the domination number of $G$. Akbari, Kumar, Mohar and Pragada \cite{akbari2024linear}, proved that $\min\{\mathcal{E}_2^+(G), \mathcal{E}_2^-(G)\}\ge \frac{3n}{4}$, whenever $n\ge 4$. But Conjecture \ref{conj:s_plus_main} still remains elusive. 

If Conjecture \ref{conj:s_plus_main} holds, then by H\"{o}lder's inequality, the lower bound of $n-1$ should naturally extend to all higher energies $\mathcal{E}^+_p(G)$ and $\mathcal{E}^-_p(G)$ for $p \geq 2$. 

In our first main result, we settle Conjecture \ref{conj:E_q-} for $p\ge 4$. 

\begin{theorem}\label{thm:negative_q_energy_main} Let $p\ge 4$. For any connected graph $G$ of order $n$ and $G\ncong K_n$, we have
\[ \mathcal{E}_p^-(G) \ge n.\]
\end{theorem}

Next, as a step toward Conjecture \ref{conj:E_q+}, we establish a lower bound for $\mathcal{E}_p^+(G)$, which is linear in $n=|V(G)|$ and scales with $p$.  

\begin{theorem}\label{thm:positive_q_energy_main}
    Let $p\ge 4$. For any connected graph $G$ of order $n \geq 4$, we have
\[ \mathcal{E}_p^+(G) \ge \bigg(\frac{4}{3}\bigg)^{p/4}n.\]
\end{theorem}

Observe that $P_n$ is bipartite and $\lambda_1(P_n) < 2$, which implies $\mathcal{E}_p^+(P_n) < 2^{p-1}n$. Hence, our lower bound has the correct order of $n$.

A key ingredient in our proofs is the following super-additivity result for the positive and negative $p$-energies of Hermitian matrices, which we prove in Section 2 using a pinching inequality for weakly unitarily invariant norms of matrices.

\begin{theorem}\label{thm:energy_matrices_main}
Let $A=[A_{ij}]_{i,j=1}^k$ be a Hermitian matrix partitioned into $k^2$ blocks, where $A_{ii}$ are square matrices. For any real number $p\geq 1$, the following hold:
$$ \mathcal{E}_p^+(A)\geq \sum_{i=1}^{k} \mathcal{E}_p^+(A_{ii}), \quad \mathcal{E}_p^-(A)\geq \sum_{i=1}^{k} \mathcal{E}_p^-(A_{ii}).$$
\end{theorem}

We note that the above result for the adjacency matrices of graphs was observed in \cite{Akbari_2009} when $p=1$, and in \cite{ akbari2024linear, Zhang_2024} when $p=2$, and was effectively used to obtain lower bounds for 1-energy and 2-energy of graphs. We believe that the above result is of independent interest and can be used to answer many questions concerning the $p$-energy of graphs.

The paper is organized as follows. We prove Theorem \ref{thm:energy_matrices_main} in Section 2. We establish Theorems \ref{thm:negative_q_energy_main} and \ref{thm:positive_q_energy_main} for $p=4$ in Section 3 and finish the proof for $p\ge 4$ in Section 4. 

\section{Super-additivity for $p$-energies}

In this section, we prove Theorem \ref{thm:energy_matrices_main}. For any $1 \leq p < \infty$, the \emph{Schatten $p$-norm} $||A\Vert_p$ of a matrix $A \in \mathbb{C}^{n\times n}$ is defined as follows: 
\[ \Vert A\Vert_p = \big[\tr(|A|^p)\big]^{1/p},\]
where $|A| = (A^*A)^{1/2}$ is the absolute value of $A$, and $\tr$ is the usual trace of matrix. 

A matrix norm $\Vert.\Vert$ in the space of matrices $\mathbb{C}^{n\times n}$ is called \textit{weakly unitarily invariant} if for all unitary matrices $U \in \mathbb{C}^{n\times n}$ and for all matrices $A \in \mathbb{C}^{n\times n}$, the norm satisfies 
\[ \Vert A\Vert = \Vert UAU^*\Vert. \]

The pinching inequality (see \cite[Chapter III]{gohberg1978introduction}, \cite[Chapter IV.4]{bhatia2013matrix}) is well-known for weakly unitarily invariant norms. For the Schatten $p$-norm, the pinching inequality takes the following form.

\begin{lemma}\label{lem:pinching}
Let $A = [A_{ij}]$ be a complex block matrix, where $1\le i,j\le k$, and let $p\ge 1$. Then 
  \[\sum_{i=1}^k \Vert A_{ii}\Vert_p^p \leq \Vert A \Vert^p_p.\]
\end{lemma}

We now proceed to the proof of Theorem \ref{thm:energy_matrices_main}.

\begin{proof}[Proof of Theorem \ref{thm:energy_matrices_main}]
Let $x_i$ denote an eigenvector corresponding to the eigenvalue $\lambda_i$ for $i=1, \ldots, n$ such that $\{x_1, \ldots, x_n\}$ forms an orthonormal basis for $\mathbb{R}^n$. By the spectral decomposition (see \cite[Theorem 4.1.5]{Horn_Johnson_2013}), we have
\[A = \sum_{i=1}^n \lambda_i x_ix_i^*.\]
Define
\[B=\sum_{\lambda_i>0}\lambda_ix_ix_i^*,\,\, C=-\sum_{\lambda_i<0} \lambda_i  x_ix_i^*.\]
Clearly, both $B$ and $C$ are positive semidefinite matrices, and the following equalities hold: 
\[BC=CB=0, \,\,\,\, A=B-C.\]
Observe that $\mathcal{E}^{+}_p(A) = \Vert B\Vert_p^p$, $\mathcal{E}^{-}_p(A) = \Vert C\Vert_p^p$.  Let $B=[B_{ij}]$ and $C=[C_{ij}], 1\leq i,j\leq k$, partitioned conformally as $A$. We have $A_{ii}=B_{ii}-C_{ii}$, for $i=1,\dots,k$. Since $B$ and $C$ are positive semidefinite matrices (see \cite[Theorem 7.7.7]{Horn_Johnson_2013}), both $B_{ii}$ and $C_{ii}$ are also positive semidefinite for $i=1,\dots,k$. By Lemma \ref{lem:pinching},  \[\Vert B \Vert^p_p \geq  \sum_{i=1}^k \Vert B_{ii}\Vert_p^p.\] 
Since $B_{ii}=A_{ii}+C_{ii}$ and $C_{ii}$ is a positive semidefinite matrix, $\lambda_r(B_{ii}) \geq \lambda_r(A_{ii})$ for $1\le r \le l$, where $l$ is the number of positive eigenvalues of $A_{ii}$. This implies 
\[\Vert B_{ii}\Vert_p^p \geq \mathcal{E}^{+}_p(A_{ii}),\] 
which yields
 \[ \mathcal{E}^{+}_p(A) = \Vert B\Vert^p_p \geq \sum_{i=1}^{k} \mathcal{E}_p^+(A_{ii}).\]
Applying the same argument to $-A$ completes the proof.
\end{proof}

\section{Positive and Negative 4-energies  of Graphs}

In this section, we establish our main results for the case $p=4$. We first investigate the positive 4-energy and show the following. 

\begin{theorem}\label{thm:positive-4-energy}
    Let $G$ be a connected graph of order $n\ge 3$. Then $\mathcal{E}_4^+(G)\ge \frac{4}{3}n$.
\end{theorem}

\begin{proof}
    We proceed by induction on $n$. The assertion has been verified using a computer for all graphs up to $8$ vertices, so assume $n\ge 9$. 

    Let $G$ be a connected graph of order $n$ with maximum degree $\Delta$. If $\Delta=2$, then $G$ is $P_n$ or $C_n$. In either case, by Theorem \ref{thm:energy_matrices_main} and induction hypothesis, we have
    \[ \mathcal{E}_4^+(G)\ge \mathcal{E}_4^+(P_3) + \mathcal{E}_4^+(P_{n-3})\geq 4 + \frac{4(n-3)}{3} = \frac{4}{3}n.\]

    So assume $\Delta \ge 3$. If $\Delta \geq \sqrt{4n/3}$, then 
    \[\mathcal{E}_4^+(G)\ge \lambda_1^4(G) \geq \Delta^2 \geq \frac{4n}{3}.\] 
    
    Hence we can assume that $3\leq \Delta < \sqrt{4n/3}$. Let $T$ be a spanning tree of $G$ rooted at a vertex $v$, where $\deg_T(v)=\Delta$. Since $n\ge 9$, we can find an edge $uv$ in $T$ such that $T-uv$ has two components $T_1$ and $T_2$, both of order at least 3. Using Theorem \ref{thm:energy_matrices_main} and the induction hypothesis,
    \[\mathcal{E}_4^+(G)\ge \mathcal{E}_4^+(G[V(T_1)]) + \mathcal{E}_4^+(G[V(T_2)])\geq \frac{4}{3}n.\] 
    This completes the proof.
 \end{proof}
 
\textbf{Remark.} With more effort one can improve the lower bound in Theorem \ref{thm:positive-4-energy}, and show that $\mathcal{E}_4^+(G)\ge 2n$ for $n\ge 5$. We omit doing this since even $2n$ is far from the optimal lower bound as given in Conjecture \ref{conj:E_q+}.

We now study the negative 4-energy of graphs. We will require the following well-known result (cf. \cite[Theorem 4.3.17]{Horn_Johnson_2013}).

\begin{theorem}[Interlacing Theorem]\label{thm:Interlacing}
Let $A$ be a real symmetric matrix of order $n$. Let $B$ be a principal submatrix of $A$ of order $n-1$. Then, for $1\le i \le n-1$,
\[\lambda_i(A) \ge \lambda_i(B) \ge \lambda_{i+1}(A).\]
\end{theorem}

We first observe the following. 

\begin{lemma}\label{lemma:disjoint_union_cliques}
Let $G\ncong K_n$ be a connected graph of order $n\ge 3$. Let $v\in V(G)$ be such that $G-v$ is a disjoint union of cliques. Then $\mathcal{E}_4^-(G)\ge n+1$.
\end{lemma}

\begin{proof}
The assertion can be verified using a computer for graphs up to 4 vertices, so assume $n\ge 5$. Assume that $G-v$ is the disjoint union of cliques $G_1, \ldots, G_\ell$. Since $K_{1,\ell}$ is an induced subgraph of $G$, $|\lambda_n(G)|\ge \sqrt{\ell}$. By the Interlacing Theorem, we have
\begin{align*}
  \mathcal{E}_4^-(G)& \ge \mathcal{E}_4^-(G-v) -\lambda_{n-1}^4(G-v) + \lambda_n^4(G) \\
  & \ge \bigg(\sum_{i=1}^{\ell} |V(G_i)|-1\bigg) - 1 + \ell^2\\
  & = n + \ell^2 -\ell -2.
\end{align*}
If $\ell\ge 3$, we immediately get $\mathcal{E}_4^-(G)\ge n+1$. We now consider the following cases:

\textbf{Case 1:} $\ell =1$. ~ Then $G$ contains an induced $P_3$ (containing $v$). By the Interlacing Theorem, 
\begin{align*}
    \mathcal{E}_4^-(G) &\ge \mathcal{E}_4^-(G-v) - \lambda_{n-1}^4(G-v) + \lambda_n^4(G)\\
    & \ge (n-2) - 1 + \lambda_{\min}^4(P_3)\\
    & = (n-3) + 4 = n+1.
\end{align*}

\textbf{Case 2:} $\ell =2$. ~ Since $n\ge 5$, $G$ contains one of the graphs $H$ in Figure \ref{fig:H} as an induced subgraph. A direct computation shows that $|\lambda_{\min}(H)|\ge 1.5$, implying $|\lambda_{\min}^4(H)|\ge 5$. As before
\begin{align*}
\mathcal{E}_4^-(G) & \ge \mathcal{E}_4^-(G-v) - \lambda_{n-1}^4(G-v) + \lambda_n^4(G)\\
& \ge (n-3) -1 + \lambda_{\min}^4(H) = n+1.
\end{align*}
The proof is complete.
\end{proof}

\begin{figure}
\begin{subfigure}{0.3\textwidth}
    \centering
\begin{tikzpicture}[scale=1.3]
\draw [line width=1pt] (-1,0)-- (-2,0.56)
 (-2,0.56)-- (-2,-0.58)
 (-2,-0.58)-- (-1,0)
 (-1,0)-- (0,0)
 (0,0)-- (1,0);

\draw [fill=black] (-1,0) circle (2pt)
 (-2,0.56) circle (2pt)
 (-2,-0.58) circle (2pt)
 (0,0) circle (2pt)
 (0,0.2) node {$v$}
 (1,0) circle (2pt);
\end{tikzpicture}
\end{subfigure}
\begin{subfigure}{0.3\textwidth}
    \centering
\begin{tikzpicture}[scale=1.3]
\draw [line width=1pt] (-1,0)-- (-2,0.56)
 (-2,0.56)-- (-2,-0.58)
 (-2,-0.58)-- (-1,0)
 (-1,0)-- (0,0.57)
 (0,0.57)-- (0,-0.59)
 (0,-0.59)-- (-1,0);
 
\draw [fill=black] (-1,0) circle (2pt)
 (-2,0.56) circle (2pt)
 (-2,-0.58) circle (2pt)
 (0,0.57) circle (2pt)
 (0,-0.59) circle (2pt)
 (-1,0.2) node {$v$};
\end{tikzpicture}
\end{subfigure}
\begin{subfigure}{0.3\textwidth}
    \centering
\begin{tikzpicture}[scale=1.3]
\draw [line width=1pt] (-1,0)-- (-2,0.56)
 (-2,0.56)-- (-2,-0.58)
 (-2,-0.58)-- (-1,0)
 (-2,0.56)-- (-3,0)
 (-3,0)-- (-2,-0.58)
 (-1,0)-- (0,0);
 
\draw [fill=black] (-1,0) circle (2pt)
 (-2,0.56) circle (2pt)
 (-2,-0.58) circle (2pt)
 (-3,0) circle (2pt)
 (0,0) circle (2pt)
 (-1,-0.2) node {$v$};
\end{tikzpicture}
\end{subfigure}\\

\begin{subfigure}{0.45\textwidth}
    \centering   
\begin{tikzpicture}[scale=1.3]
\draw [line width=1pt] (-2,0)--(-1,0)--(0,0)--(1,0)--(2,0);

\draw [fill=black] (-2,0) circle (2pt)
 (-1,0) circle (2pt)
 (0,0) circle (2pt)
 (1,0) circle (2pt)
 (2,0) circle (2pt)
 (0,-0.2) node {$v$};
\draw [color=white] (0,-1) circle (2pt);
\end{tikzpicture}

\end{subfigure}
\begin{subfigure}{0.45\textwidth}
    \centering
\begin{tikzpicture}[scale=1.3]
\draw [line width=1pt] (-2,0)-- (-1,1)
 (-1,1)-- (0,0)
 (0,0)-- (-1,-1)
 (-1,-1)-- (-2,0)
 (-2,0)-- (0,0)
 (-1,1)-- (-1,-1)
 (0,0)-- (1.39,0);

\draw[fill=black] (-2,0) circle (2pt)
 (-1,1) circle (2pt)
 (0,0) circle (2pt)
 (-1,-1) circle (2pt)
 (1.39,0) circle (2pt)
 (0,-0.2) node {$v$};
\end{tikzpicture}
\end{subfigure}
    \caption{Induced subgraph(s) $H$}
    \label{fig:H}
\end{figure}

Next, we give a lower bound on the negative $4$-energy of graphs with a dominating vertex. Recall that a vertex is called \emph{dominating} if it is adjacent to every other vertex in the graph. 

\begin{theorem}\label{thm:4_energy_dominating}
Let $G \ncong K_n$ be a connected graph of order $n\geq 3$ with a dominating vertex $v$. Then $\mathcal{E}_4^-(G) \geq n+1$.
\end{theorem}

\begin{proof}
If $n=3$, then $G \cong P_3$, and $\mathcal{E}_4^-(P_3) = 4$. If $n=4$, one can check manually (or using a computer) that the stronger claim $\mathcal{E}_4^-(G)\ge 6$ is true. So, assume $n\ge 5$. We proceed by induction on $n$. 

If $G-v$ is a disjoint union of cliques, then we are done by Lemma \ref{lemma:disjoint_union_cliques}. So, assume there is an induced path $abc$ in $G-v$. If $G-\{a,b,c\}$ is not a complete graph, then by Theorem \ref{thm:energy_matrices_main} and induction hypothesis, we have 
\[\mathcal{E}_4^-(G) \geq \mathcal{E}_4^-(P_3) + \mathcal{E}_4^-(G-\{a,b,c\}) \ge 4 + (n-2)> n+1.\]
So suppose $G-\{a,b,c\}$ is complete. Then $G[a,b,c,v]$ has a dominating vertex $v$ and $G-\{a,b,c,v\}\cong K_{n-4}$. By Theorem \ref{thm:energy_matrices_main}, we have 
 \[\mathcal{E}_4^-(G) \geq \mathcal{E}_4^-(G[a,b,c,v]) + \mathcal{E}_4^-(K_{n-4}) \ge 6 + (n-5)\ge n+1.\] 
The proof is complete.
\end{proof}

We are now ready to prove Theorem \ref{thm:negative_q_energy_main} when $p=4$.

\begin{theorem}\label{thm:negative-4-energy}
Let $G\ncong K_n$ be a connected graph of order $n\ge 3$. Then $\mathcal{E}_4^-(G)\ge n$.
\end{theorem}

\begin{proof}
We proceed by induction on $n$. The assertion can be verified for all graphs up to 8 vertices, so assume $n\ge 9$. Throughout the proof, we use Theorem \ref{thm:energy_matrices_main} and the induction hypothesis repeatedly, without mention. 

Let $G\ncong K_n$ be a connected graph of order $n$ with maximum degree $\Delta$. If $\Delta=2$, then $G$ is $P_n$ or $C_n$. In either case
\[ \mathcal{E}_4^-(G)\ge \mathcal{E}_4^-(P_3) + \mathcal{E}_4^-(P_{n-3})\ge n.\]

Henceforth, assume that $\Delta\ge 3$. Let $v\in V(G)$ be such that $\deg(v)=\Delta$. Suppose $C$ is a component of $G-N[v]$ such that $C$ is non-complete. Clearly, $C$ has at least 3 vertices. Also, $G-V(C)$ is connected and has order at least 3, since $\Delta \ge 3$. Moreover, if $G-V(C)$ is complete, then $G-V(C)\cong K_t$, where $t \geq \Delta+1$, and so $G-V(C)$ and $C$ are vertex disjoint, a contradiction. So $G-V(C)$ is non-complete. We have
\[\mathcal{E}_4^-(G) \ge \mathcal{E}_4^-(C) + \mathcal{E}_4^-(G-C)\ge n. \]

So we can assume that $G-N[v]$ is a disjoint union of cliques $C_1, \ldots,C_k$ for some $k\ge 0$. If $k=0$, then $v$ is a dominating vertex in $G$, and we are done by Theorem \ref{thm:4_energy_dominating}. So assume $k\ge 1$. We consider the following cases:

\textbf{Case 1:} There exists $u\in N(v)$ such that $u$ has neighbours in at least two components of $G-N[v]$.

Without loss of generality, let $C_1, \ldots, C_t$ be cliques of $G-N[v]$ such that $u$ has at least one neighbour in $C_i$ for $1\le i\le t$. Since $t\ge 2$, $G[V(C_1\cup \cdots \cup C_t)\cup \{u\}]$ is a non-complete graph. If $G-V(C_1\cup \cdots \cup C_t) -u$ is non-complete, then
\[\mathcal{E}_4^-(G) \ge \mathcal{E}_4^-(G[V(C_1\cup \cdots \cup C_t)\cup \{u\}]) + \mathcal{E}_4^-(G-V(C_1\cup \cdots \cup C_t) -u) \ge n. \]
If $G-V(C_1\cup \cdots \cup C_t) - u$ is a complete graph, then $G-u$ is a disjoint union of cliques and we are done by Lemma \ref{lemma:disjoint_union_cliques}. 

\textbf{Case 2:} For any vertex $u\in N(v)$, there is at most one $C_i$ such that $N(u)\cap V(C_i)\ne \emptyset$.

Let $u\in N(v)$ be such that $u$ has a neighbour in a unique $C_i$. Suppose $G[V(C_i)\cup \{u\}]$ is non-complete. If $G-V(C_i) - u$ is non-complete, then 
\[\mathcal{E}_4^-(G) \ge \mathcal{E}_4^-(G[V(C_i)\cup \{u\}]) + \mathcal{E}_4^-(G-V(C_i)-u)\ge n. \]
If $G-V(C_i)-u$ is a complete graph, then we are done using Lemma \ref{lemma:disjoint_union_cliques}, since $G-u$ is a disjoint union of cliques. So we can assume that $G[V(C_i)\cup \{u\}]$ is complete, whenever a vertex $u\in N(v)$ has a neighbour in $C_i$, for some $i$. 

Define $S_i=\{u\in N(v): u \text{ has a neighbour in }C_i\}$. Now, suppose $G[V(C_i)\cup S_i]$ is non-complete. Since every vertex $w\in V(C_i)$ is a dominating vertex of $G[V(C_i)\cup S_i]$, using Theorem \ref{thm:4_energy_dominating} we have
\[\mathcal{E}_4^-(G) \ge \mathcal{E}_4^-(G[V(C_i)\cup S_i]) + \mathcal{E}_4^-(G-(V(C_i)\cup S_i))\ge n.\] Hence, to proceed, we assume that $G[V(C_i)\cup S_i]$ is complete for all $i$. 

Now, if $G-v$ is a disjoint union of cliques, then we can apply Lemma \ref{lemma:disjoint_union_cliques} to argue that $\mathcal{E}_4^-(G)\ge n$. So assume that $G-v$ is not a disjoint union of cliques. Using the observations above, this can happen only if one of the following occurs:

\textbf{Subcase 2.1:} There exist vertices $a\in S_i$ and $b\in S_j$, where $i\ne j$, such that $a$ is adjacent to $b$. Pick $a'\in V(C_i)$ and $b'\in V(C_j)$. Then $a'abb'$ is an induced $P_4$ in $G$. We have
\begin{align*}
 \mathcal{E}_4^-(G) & \ge \mathcal{E}_4^-(P_4) + \mathcal{E}_4^-(C_i - a') + \mathcal{E}_4^-(C_j - b') + \mathcal{E}_4^-(G-V(C_i\cup C_j)-a-b)\\
 & \ge 7 + (|V(C_i)|-2) + (|V(C_j)|-2) + (n - |V(C_j)|-|V(C_j)|-3)\\
 &= n.
\end{align*}

\textbf{Subase 2.2:} There exists a vertex $a\in S_i$ and a vertex $b\in N(v)$ such that $b$ has no neighbours in $G-N[v]$ and $a$ is adjacent to $b$. Then $a$ is a dominating vertex in the non-complete graph $G[V(C_i)\cup S_i\cup \{b\})]$. Using Theorem \ref{thm:4_energy_dominating}, we have 
\begin{align*}
  \mathcal{E}_4^-(G) &\ge \mathcal{E}_4^-(G[V(C_i)\cup S_i\cup \{b\})]) + \mathcal{E}_4^-(G-(V(C_i)\cup S_i\cup\{b\}))\\
  & \ge (|V(C_i)|+|S_i|+2) + (n-|V(C_i)|-|S_i|-2) =n.  
\end{align*}

\textbf{Subase 2.3:} There exist vertices $a,b,c\in N(v)$ such that $abc$ is an induced path in $G$ and none of the vertices $a,b,c$ have neighbours in $G-N[v]$. ~ We have 
\[\mathcal{E}_4^-(G) \ge \mathcal{E}_4^-(G[a,b,c]) + \mathcal{E}_4^-(G-a-b-c))\ge 4 + (n-4) = n.\]
This completes the proof.
\end{proof}

\section{Proof of Theorems \ref{thm:negative_q_energy_main} and \ref{thm:positive_q_energy_main}}

 We first recall a consequence of H\"{o}lder's inequality. For any $x = (x_1,\dots,x_n) \in \mathbb{R}^n$ and any $1 \leq p < \infty$, the $p$-norm of $x$ is given by
 \[\Vert x \Vert_p = \Big(\sum_{i=1}^{n} \vert x_i \vert^p \Big)^\frac{1}{p}.\]
 The following lemma provides a fundamental inequality relating different $p$-norms. 
 
\begin{lemma}\label{lem:p_norm_ineq}
Let $x$ be a vector in $\mathbb{R}^n$ and $1 \leq q\leq p$. Then
 \[ \Vert x\Vert_p \leq \Vert x\Vert_q \leq n^{\frac{1}{q} -\frac{1}{p}} \Vert x\Vert_p.\]
\end{lemma}

The next lemma follows easily from Lemma \ref{lem:p_norm_ineq} and compares $p$-energy values for different $p$.

\begin{lemma}\label{scaling_lemma}
    Let $A$ be an Hermitian matrix and let $n^+,n^-$ denote the number of positive and negative eigenvalues of $A$, respectively. Then for any real numbers $p\geq q\geq 1$, we have
   \begin{align*}
           \frac{\Big(\mathcal{E}_q^+(A)\Big)^{p/q}}{(n^+)^{\frac{p}{q}-1}}\leq \mathcal{E}_p^+(A) \leq \Big(\mathcal{E}_q^+(A)\Big)^{p/q}.
   \end{align*}
   The same inequalities hold for $\mathcal{E}_p^-(A)$ with $n^-$.
\end{lemma}

We are now ready to prove Theorems \ref{thm:negative_q_energy_main} and \ref{thm:positive_q_energy_main}.

\begin{proof}[Proof of Theorems \ref{thm:negative_q_energy_main} and \ref{thm:positive_q_energy_main}]
    By the first inequality of Lemma \ref{scaling_lemma} with $p=4$ and Theorem \ref{thm:positive-4-energy} and noting that $n^+ \leq n$, we get 
    \begin{align*}
        \mathcal{E}_p^+(G) &\geq \frac{\Big(\mathcal{E}_4^+(G)\Big)^{\frac{p}{4}}}{(n^+)^{\frac{p}{4}-1}} \\
        &\ge \frac{(\frac{4}{3}n)^{\frac{p}{4}}}{(n^+)^{\frac{p}{4}-1}} \ge \bigg(\frac{4}{3}\bigg)^{\frac{p}{4}}n.
    \end{align*}
    This completes the proof for $\mathcal{E}_p^+(G)$. The proof for $\mathcal{E}_p^-(G)$ is similar. 
\end{proof}

\textbf{Remark.} \label{remark:upper_bound}
 We mention the following upper bounds for the positive and negative 2-energy of graphs, which in turn yield upper bounds for $\mathcal{E}_p^+(G)$ and $\mathcal{E}_p^-(G)$.
 
\begin{proposition}[\cite{Elphick_FGW_2016, Elphick_Linz_2024}]\label{prop:Elphick_Sminus_bnd}
    Let $G$ be a graph of order $n$. Then 
    \[\mathcal{E}^+_2(G) \leq (n-1)^2, \quad \mathcal{E}^-_2(G) \leq \Big(\frac{n}{2}\Big)^2.\]
\end{proposition}

Using the inequalities from the Proposition \ref{prop:Elphick_Sminus_bnd} along with the second inequality from Lemma \ref{scaling_lemma} for $p=2$, we get
\[ \mathcal{E}_p^+(G)\le (n-1)^p, \quad \mathcal{E}^-_p(G) \leq \Big(\frac{n}{2}\Big)^p.\]

\section*{Acknowledgements}
The authors would like to thank Quanyu Tang for his fruitful comments in the preparation of the paper. The authors also thank the anonymous referee for their valuable feedback. 

\bibliographystyle{plainurl}
\bibliography{p_references}

\vspace{0.6cm}
\noindent Saieed Akbari, Email: {\tt s\_akbari@sharif.edu}\\
The research visit of Saieed Akbari at Simon Fraser University was supported in part by the ERC Synergy grant (European Union, ERC, KARST, project number 101071836).\\ 
\textsc{Dept. of Mathematical Sciences, Sharif University of Technology, Tehran, Iran}\\[2pt]

\noindent Hitesh Kumar, Email: {\tt hitesh.kumar.math@gmail.com}, {\tt hitesh\_kumar@sfu.ca}\\
\textsc{Dept. of Mathematics, Simon Fraser University, Burnaby, BC \ V5A 1S6, Canada}\\[2pt]

\noindent Bojan Mohar, Email: {\tt mohar@sfu.ca}\\
Bojan Mohar is supported in part by the NSERC Discovery Grant R832714 (Canada), by the ERC Synergy grant (European Union, ERC, KARST, project number 101071836), and by the Research Project N1-0218 of ARIS (Slovenia). On leave from FMF, Department of Mathematics, University of Ljubljana.\\
\textsc{Dept. of Mathematics, Simon Fraser University, Burnaby, BC \ V5A 1S6, Canada}\\[2pt]

\noindent Shivaramakrishna Pragada, Email: {\tt shivaramakrishna\_pragada@sfu.ca}\\
\textsc{Dept. of Mathematics, Simon Fraser University, Burnaby, BC \ V5A 1S6, Canada}

\end{document}